\documentclass[12pt]{article}

\usepackage[T1]{fontenc}
\usepackage[english]{babel}
\usepackage[a4paper,margin=1in]{geometry}
\usepackage{amsmath,amssymb,amsthm,mathtools,bm}
\usepackage{booktabs,array}
\usepackage{microtype}
\usepackage{hyperref}
\usepackage{enumitem}
\usepackage{cleveref}
\usepackage{caption}
\usepackage{float}
\usepackage{tikz}
\usetikzlibrary{arrows.meta,positioning,fit,calc,decorations.pathreplacing}
\allowdisplaybreaks
\hypersetup{
	hidelinks,
	pdftitle={Counterexamples and Optimal Order Bounds for Free Copies under Edge Set Mappings},
	pdfsubject={Set mappings, free copies, probabilistic embeddings, and extremal graphs},
	pdfkeywords={set mapping, free copy, probabilistic method, Lovasz local lemma, extremal graph theory}
}

\newtheorem{theorem}{Theorem}[section]
\newtheorem{lemma}[theorem]{Lemma}
\newtheorem{conjecture}[theorem]{Conjecture}
\newtheorem{corollary}[theorem]{Corollary}
\newtheorem{proposition}[theorem]{Proposition}

\theoremstyle{remark}

\newcommand{\PP}{\mathbb{P}}

\newcommand{\calF}{\mathcal{F}}

\crefname{theorem}{Theorem}{Theorems}
\crefname{lemma}{Lemma}{Lemmas}
\crefname{corollary}{Corollary}{Corollaries}
\crefname{conjecture}{Conjecture}{Conjectures}
\crefname{remark}{Remark}{Remarks}
\crefname{section}{Section}{Sections}

\title{Sharp Bounds for Edge Mappings of Graphs}

\author{
	Zhanping Yang\thanks{Corresponding author. Center for Discrete Mathematics, Fuzhou University, Fuzhou, Fujian 350108, China. Research supported by National Key R\&D Program of China (Grant No. 2023YFA1010202) and National Natural Science Foundation of China (Grant No. 12371342). Email: \texttt{yangzp@fzu.edu.cn}.}}
\date{}

\begin{document}
	\maketitle
	
	\begin{abstract}
The study of extremal problems for set mappings was initiated by Erd\H{o}s and Hajnal in 1958. Alon and Caro, and later Caro, Patk\'os, Tuza and Vizer, studied related problems for edge mappings of graphs. For a graph $G$ with $n$ vertices and $m$ edges, let $g_0(G)$ be the minimum integer $N$ such that every mapping $f:E(K_N)\to E(K_N)$ satisfying $f(e)\cap e=\varnothing$ for every $e\in E(K_N)$ admits a copy $G^*$ of $G$ in $K_N$ with $f(e)\notin E(G^*)$ for every $e\in E(G^*)$. Gishboliner, Jin and Sudakov conjectured that $g_0(G)=O(n+m/\sqrt n)$. In this paper, we disprove this conjecture by constructing a family of graphs. We then prove that $g_0(G)=O(n+m^{3/4})$. Furthermore, for every $n\ge2$ and $m$ satisfying $\lceil n/2 \rceil\leq m\leq \binom{n}{2}$, there exists a graph $G_*$ with $n$ vertices, $m$ edges and no isolated vertices satisfying $g_0(G_*)=\Omega(n+m^{3/4})$. Thus the upper bound is best possible in terms of $n$ and $m$.
	\end{abstract}
	
	
	
	\section{Introduction}\label{sec:intro}
Set mapping problems were introduced by Erd\H{o}s and Hajnal~\cite{ErdosHajnal1958}. For fixed positive integers $k,\ell$ and an integer $N$, consider all mappings $f:\binom{[N]}{k}\longrightarrow\binom{[N]}{\ell}$ satisfying $f(X)\cap X=\varnothing$ for every
$X\in\binom{[N]}{k}$. Let $p(f)$ be the maximum size of a set $P\subseteq[N]$ such that $f(X)\cap P=\varnothing$ for every $X\in\binom{P}{k}$, and define $p_{k,\ell}(N)=\min_f p(f)$, where the minimum is taken over all such mappings. Erd\H{o}s and Hajnal~\cite{ErdosHajnal1958} proved that $\Omega\bigl(N^{1/(k+1)}\bigr)\le p_{k,\ell}(N)\le O\bigl((N\log N)^{1/k}\bigr)$. Spencer~\cite{Spencer1972} improved the lower bound to $\Omega(N^{1/k})$. F\"uredi~\cite{Furedi1991} proved that $p_{2,\ell}(N)=\Theta(N^{1/2})$. Later, Conlon, Fox and Sudakov~\cite{CFS2016} proved that $p_{k,\ell}(N)=\Theta(N^{1/k})$ whenever $\ell\ge(k-1)!$. Alon, Caro and Tuza~\cite{AlonCaroTuza1989} also studied related problems concerning sub-Ramsey numbers.
		
In the graph setting, Alon and Caro~\cite{AlonCaro1986} studied edge mappings $f:E(K_N)\to E(K_N)$. A subgraph $H$ is called $f$-fixed if $f(e)=e$ for every $e\in E(H)$, and $f$-free if $f(e)\notin E(H)$ for every $e\in E(H)$. They proved that for every $r,t\ge1$, the largest $N$ for which there exists an edge mapping of $K_N$ with neither an $f$-fixed $K_{1,r}$ nor an $f$-free $K_{1,t}$ is $2r+t-2$. Caro~\cite{Caro1987} extended these problems to $k$-uniform hypergraphs and proved that $g(s,k,0)=s^{2-1/k+o(1)}$, where $g(s,k,0)$ is the largest $N$ for which there exists a mapping $f:\binom{[N]}{k}\to\binom{[N]}{k}$ satisfying $f(e)\cap e=\varnothing$ for every $e\in\binom{[N]}{k}$, with no $f$-free copy of $K_s^{(k)}$. Here $K_s^{(k)}$ denotes the complete $k$-uniform hypergraph on $s$ vertices. Caro and Sch\"onheim~\cite{CaroSchonheim1988}
also investigated edge mappings with fixed edges under the condition that no $f$-free copy of $K_3$ exists.

More recently, Caro, Patk\'os, Tuza and Vizer~\cite{CPTV2026} studied the existence of $f$-free copies of general graphs. For $a\in\{0,1\}$, let $\mathcal F_{N,a}$ denote the family of mappings
$f:E(K_N)\to E(K_N)$ satisfying $|f(e)\cap e|\le a$ for every $e\in E(K_N)$. For a graph $G$ with $n$ vertices, let $g_a(G)$ be the minimum integer $N\ge\max\{n,4\}$ such that every $f\in\mathcal F_{N,a}$ admits an $f$-free copy of $G$ in $K_N$. They also considered the parameter $w(G)$, which requires a copy $G^*$ satisfying $f(e)\cap V(G^*)=\varnothing$ for every $e\in E(G^*)$. For every graph $G$ with $n\ge4$ vertices and $m\ge1$ edges, they proved $g_0(G)\le w(G)\le2nm-4m+2$. Moreover, they determined $g_0(tK_2)=g_1(tK_2)=2t$ for $t\ge4$. In a related paper, the same authors~\cite{CPTV2025} studied $f$-free matchings in general host graphs and proved that, for every $t\ge3$ and sufficiently large $N$, the maximum number of edges in an $N$-vertex graph admitting an edge mapping $f$ with $f(e)\ne e$ for every edge $e$ and no $f$-free copy of $tK_2$ is $\operatorname{ex}(N,tK_2)+t-1$.

Gishboliner, Jin and Sudakov~\cite{GJS2026} improved the upper bounds on $w(G)$. For every graph $G$ with $n$ vertices, $m$ edges and no isolated vertices, they proved $w(G)=O(m)$, $g_1(G)=O(m)$ and 
\begin{equation*}
g_0(G)=\Omega\left(n+\frac{m}{\sqrt{n\log n}}\right).
\end{equation*}
In particular, $g_0(G)\le w(G)$ gives $g_0(G)=O(m)$. For complete graphs, the results of Conlon, Fox and Sudakov~\cite{CFS2016} imply $g_0(K_n)=\Theta(n^{3/2})$ and $w(K_n)=\Theta(n^2)$. Motivated by the lower bound above and the complete graph case, Gishboliner, Jin and Sudakov~\cite{GJS2026} proposed the following conjecture.
\begin{conjecture}[\textbf{Gishboliner--Jin--Sudakov}~\cite{GJS2026}]\label{conj:original}
	Let \(G\) be a graph with \(n\) vertices, \(m\) edges and no isolated vertices. Then 
	\begin{equation*}
		g_0(G)=O\!\left(n+\frac{m}{\sqrt n}\right).
	\end{equation*}
\end{conjecture}
Using the Lov\'asz local lemma, Gishboliner, Jin and Sudakov~\cite{GJS2026} proved that $g_0(G)=O\left(n+\sqrt{\Delta(G)m}\right)$, which confirms the conjecture when $\Delta(G)=O(m/n)$. 

In this paper, we disprove Conjecture~\ref{conj:original} by constructing a family of graphs $G_t$ satisfying $|V(G_t)|=\Theta(t^8)$, $|E(G_t)|=\Theta(t^{12})$ and $g_0(G_t)>t^9$. Consequently,
\begin{equation*}
\frac{g_0(G_t)}{|V(G_t)|+|E(G_t)|/\sqrt{|V(G_t)|}}\longrightarrow\infty.
\end{equation*}
We then establish a general upper bound and show that its order is best possible among bounds depending only on $n$ and $m$.
	\begin{theorem}\label{thm:upper}
		Let \(G\) be a graph with \(n\) vertices, \(m\) edges and no isolated vertices. Then 
		\begin{equation*}
			g_0(G)=O\left(n+m^{3/4}\right).
		\end{equation*}
	\end{theorem}
	The following result shows that the bound in Theorem~\ref{thm:upper} is best possible for every admissible pair of $n$ and $m$.
	\begin{theorem}\label{thm:lower-all}
		For every $n\ge2$ and every integer $m$ satisfying $\lceil n/2\rceil\le m\le\binom n2$, there exists a graph $G_*$ with $n$ vertices, $m$ edges
		and no isolated vertices such that
		\begin{equation*}
			g_0(G_*)=\Omega\left(n+m^{3/4}\right).
		\end{equation*}
	\end{theorem}
The paper is organized as follows. Section~\ref{sec:Auxiliary Lemmas and Counterexamples} establishes preliminary results and constructs a family of graphs disproving Conjecture~\ref{conj:original}. Section~\ref{sec:The Upper Bound} proves the upper bound in Theorem~\ref{thm:upper}. Section~\ref{sec:The Lower Bound} establishes the lower bound in Theorem~\ref{thm:lower-all}.
	\section{Auxiliary Lemmas and Counterexamples}\label{sec:Auxiliary Lemmas and Counterexamples}
We first show that $g_0(H)\le g_0(G)$ whenever $H$ is a subgraph of $G$, and give a form of the Lov\'asz local lemma. We then apply a result of Conlon, Fox and Sudakov~\cite{CFS2016} on complete graphs to construct a family of connected graphs disproving Conjecture~\ref{conj:original}.
	\begin{proposition}\label{prop:mono}
		If $H$ is a subgraph of $G$, then $g_0(H)\le g_0(G)$.
	\end{proposition}
	\begin{proof}[\bf Proof]
		Let $N=g_0(G)$ and fix any $f\in\mathcal F_{N,0}$. By the definition of $g_0(G)$, there exists an $f$-free copy $G^*$ of $G$ in $K_N$. Since $H\subseteq G$, there is a copy $H^*$ of $H$
		with $H^*\subseteq G^*$. Hence
		\begin{equation*}
			f(E(H^*))\cap E(H^*)\subseteq f(E(G^*))\cap E(G^*)=\varnothing.
		\end{equation*}
		Thus $H^*$ is $f$-free, and $g_0(H)\le N=g_0(G)$.
	\end{proof}
	To establish the upper bound on $g_0(G)$, we also use the following form of the Lov\'asz local lemma.
	\begin{lemma}\label{lem:variable-lll}
		Let $X_1,\ldots,X_s$ be independent random variables, and let $\mathcal B$ be a finite family of events. Suppose that each $B\in\mathcal B$ depends only on $\{X_j:j\in I_B\}$, where $I_B\subseteq[s]$ and $1\leq|I_B|\le4$. If for every $j\in[s]$,
		\begin{equation}\label{eq:PB<=1/16}
			\sum_{\substack{B\in\mathcal B,\ j\in I_B}}\mathbb P(B)\le\frac1{16},
		\end{equation}
		then
		\begin{equation*}
			\mathbb P\left(\bigcap_{B\in\mathcal B}\overline B\right)>0.
		\end{equation*}
	\end{lemma}
	\begin{proof}[\bf Proof]
	For distinct events $B,C\in\mathcal B$, write $C\sim B$ if $I_B\cap I_C\ne\varnothing$. Since the random variables are independent, each event $B$ is independent of all events determined by variables outside $I_B$. For each $B\in\mathcal B$, set $x_B=2\mathbb P(B)$. By \eqref{eq:PB<=1/16}, we have $0\le x_B\le1/8$. Moreover,
		\begin{equation*}
			\sum_{C\sim B}x_C=2\sum_{C\sim B}\mathbb P(C)\le2\sum_{j\in I_B}\sum_{\substack{C\in\mathcal B,\ j\in I_C}}\mathbb P(C)\le\frac{|I_B|}{8}\le\frac12.
		\end{equation*}
		Consequently,
		\begin{equation*}
			x_B\prod_{C\sim B}(1-x_C)\ge x_B\left(1-\sum_{C\sim B}x_C\right)\ge\frac{x_B}{2}=\mathbb P(B).
		\end{equation*}
		By Lov\'asz local lemma~\cite{AlonSpencer2016}, we obtain
		\begin{equation*}
			\mathbb P\left(\bigcap_{B\in\mathcal B}\overline B\right)>0,
		\end{equation*}
		as required.
	\end{proof}
The construction of Conlon, Fox and Sudakov~\cite[Theorem~2.2(ii)]{CFS2016} gives the following lower bound for complete graphs, which will be used to disprove Conjecture~\ref{conj:original}.
	\begin{theorem}[\textbf{Conlon–Fox–Sudakov}~\cite{CFS2016}]\label{thm:grid}
		For every integer $q\ge2$, we have 
		\begin{equation*}
			g_0\bigl(K_{3q^2+1}\bigr)>q^3.
		\end{equation*}
	\end{theorem}
	The following lower bound for complete graphs follows from Theorem~\ref{thm:grid}.
	\begin{corollary}\label{cor:clique}
		For every integer $s\ge4$, we have
		\begin{equation*}
			g_0(K_s)\ge\frac{s^{3/2}}{200}.
		\end{equation*}
	\end{corollary}
	\begin{proof}[\bf Proof]
		For $4\le s\le12$, the result follows from $g_0(K_s)\ge s\ge{s^{3/2}}/{200}$. Assume that $s\ge13$ and let $q=\lfloor\sqrt{{(s-1)}/{3}}\rfloor$. Then $q\ge2$ and $3q^2+1\le s$. By Proposition~\ref{prop:mono} and Theorem~\ref{thm:grid}, we obtain
		\begin{equation*}
			g_0(K_s)\ge g_0(K_{3q^2+1})>q^3.
		\end{equation*}
		Since $\sqrt{(s-1)/3}\ge2$, we have $2q\ge\sqrt{{(s-1)}/{3}}$. Hence
		\begin{equation*}
			q^3\ge\frac18\left(\frac{s-1}{3}\right)^{3/2}\ge \frac{s^{3/2}}{8\cdot6^{3/2}}>\frac{s^{3/2}}{200}.
		\end{equation*}
		This completes the proof.
	\end{proof}

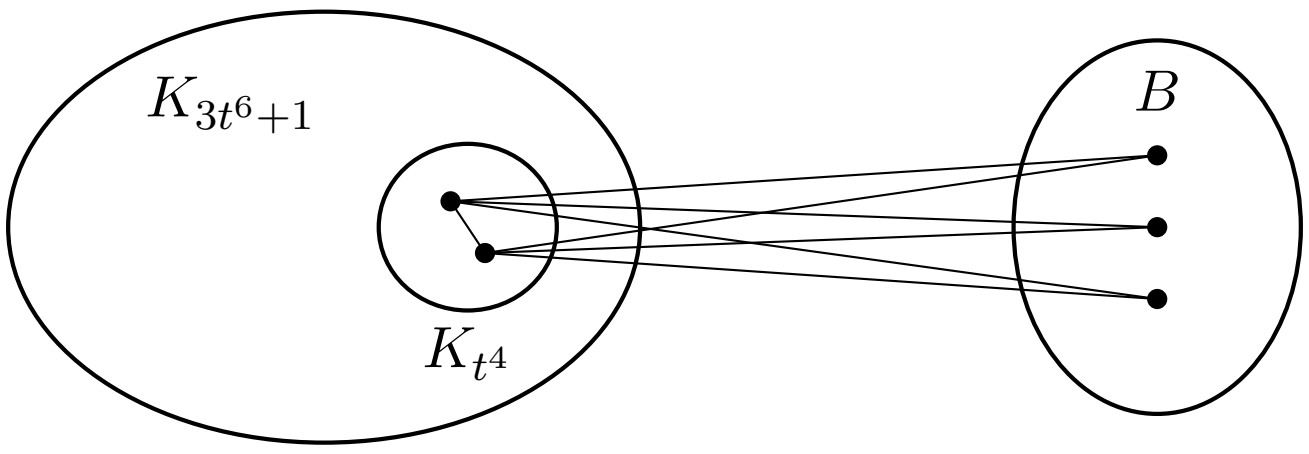
\begin{figure}[H]
		\centering
		\begin{tikzpicture}[scale=0.95, every node/.style={font=\small}]
			
			\draw[thick] (0,0) ellipse (2.2 and 1.5);
			\node at (-0.65,0.85) {$K_{3t^6+1}$};
			
			\draw[thick] (1.0,0) ellipse (0.62 and 0.58);
			\node at (1.0,-0.88) {$K_{t^4}$};
			
			\draw[thick] (5.8,0) ellipse (1.0 and 1.3);
			\node at (5.8,0.95) {$B$};
			
			\fill (0.88,0.18) circle (2pt);
			\fill (1.12,-0.18) circle (2pt);
			
			\draw (0.88,0.18) -- (1.12,-0.18);
			
			\foreach \y in {0.5,0,-0.5} {
				\fill (5.8,\y) circle (2pt);}
			
			\foreach \ax/\ay in {0.88/0.18,1.12/-0.18} {
				\foreach \by in {0.5,0,-0.5} {
					\draw (\ax,\ay) -- (5.8,\by);}}
		\end{tikzpicture}
		\caption{The graph $G_t$. }
		\label{fig:gt-construction}
	\end{figure}

	For every integer $t \geq 2$, the connected graph $G_t$ is constructed as follows (see Figure~\ref{fig:gt-construction}). Take a complete graph $K_{3t^6+1}$ and choose a set $A \subseteq V(K_{3t^6+1})$ with $|A| = t^4$. Add an independent set $B$ of size $t^8$ and join every vertex of $A$ to every vertex of $B$, adding no other edges.
	\noindent Thus $G_t$ satisfies $|V(G_t)|=t^8+3t^6+1$ and
	\begin{equation*} |E(G_t)|=\binom{3t^6+1}{2}+|A||B|=\frac{(3t^6+1)3t^6}{2}+t^4t^8=\frac{11}{2}t^{12}+\frac32t^6.
	\end{equation*}
	\begin{proposition}\label{prop:connected-counterexample}
		For every integer $t\ge2$, the graph $G_t$ satisfies
		\begin{equation*}
			\frac{g_0(G_t)}{|V(G_t)|+|E(G_t)|/\sqrt{|V(G_t)|}}=\Omega(t)\longrightarrow\infty.
		\end{equation*}
	\end{proposition}
	\begin{proof}[\bf Proof]
		Since $K_{3t^6+1}\subseteq G_t$, Proposition~\ref{prop:mono} and~Theorem~\ref{thm:grid}, applied with $q=t^3$, give
		\begin{equation*}
			g_0(G_t)\ge g_0(K_{3t^6+1})>t^9.
		\end{equation*}
		On the other hand, $|V(G_t)|=\Theta(t^8)$ and $|E(G_t)|=\Theta(t^{12})$, so
		\begin{equation*}
			|V(G_t)|+\frac{|E(G_t)|}{\sqrt{|V(G_t)|}}=\Theta(t^8).
		\end{equation*}
		It follows that
		\begin{equation*}
			\frac{g_0(G_t)}{|V(G_t)|+|E(G_t)|/\sqrt{|V(G_t)|}}=\Omega(t)\longrightarrow\infty.
		\end{equation*}
		This completes the proof.
	\end{proof}
Consequently, Conjecture~\ref{conj:original} fails even for connected graphs.
	
	\section{The Upper Bound}\label{sec:The Upper Bound}
	In this section, we prove Theorem~\ref{thm:upper}. We first find disjoint sets $S,R\subseteq V(K_N)$ such that $K_N[S]$ is $f$-free and the edge mapping satisfies certain bounds on $S$.
	We then embed the vertices of large degree into $S$ and the remaining vertices randomly into $R$, and apply Lemma~\ref{lem:variable-lll}.
	
	Fix $f\in\mathcal F_{N,0}$. For each $x\in[N]$, define
	\begin{equation}\label{eq:af}
		a_f(x)=\bigl|\{e\in E(K_N):x\in f(e)\}\bigr|.
	\end{equation}
	For distinct $x,y\in[N]$, define
	\begin{equation}\label{eq:bf}
		b_f(x,y)
		=\bigl|\{z\in[N]\setminus\{x\}: y\in f(xz)\}\bigr|.
	\end{equation}
	
	\begin{lemma}\label{lem:host-core}
		Let $q\ge1$ be an integer and suppose that  $N\ge256q^{3/2}$. For every $f\in\calF_{N,0}$, there exist disjoint sets $S,R\subseteq[N]$ with $|S|=q$ and $|R|\ge N/2$ satisfying the following properties:
		\begin{enumerate}[label=(\alph*)]
			\item  $K_N[S]$ is $f$-free;\label{a}
			\item  For every $e\in E(K_N)$, if $|(e\cup f(e))\cap S|=3$, then $(e\cup f(e))\cap R=\varnothing$;\label{b}
			\item  For every $x\in S$, $a_f(x)\le2N$;\label{c}
			\item  $\bigl|\{e\in E(K_N):f(e)\subseteq S\}\bigr|\le512q^2$;\label{d}
			\item  For every $x\in S$,\label{e}
			\begin{equation}\label{eq:H5}
				\sum_{y\in S\setminus\{x\}}\bigl(b_f(x,y)+b_f(y,x)\bigr)
				\le4096q.
			\end{equation}
		\end{enumerate}
	\end{lemma}
	\begin{proof}[\bf Proof]
		Since every $f(e)$ contains two vertices, double counting gives
		\begin{equation*}
			\sum_{x\in[N]}a_f(x)=2\binom N2=N(N-1)<N^2.
		\end{equation*}
		Hence, the set $X:=\{x\in[N]:a_f(x)\le2N\}$ satisfies $|X|>N/2$. Let $L=4q$. Since $|X|>N/2\ge128q^{3/2}\ge4q=L$. We can choose an $L$-element subset $Y\subseteq X$ uniformly at random. 

For any fixed subset $V'\subseteq[N]$ of $k$ vertices, we have $\mathbb P(V'\subseteq Y)=0$ if
$V'\nsubseteq X$ or $k>L$. Otherwise, $V'\subseteq X$ and $k\le L$, so
		\begin{equation}\label{eq:subsetprob}
			\mathbb{P}(V' \subseteq Y)=\frac{\binom{|X|-k}{L-k}}{\binom{|X|}{L}} =\prod_{i=0}^{k-1} \frac{L - i}{|X| - i} \leq \left(\frac{L}{|X|}\right)^k \leq \left(\frac{2L}{N}\right)^k.
		\end{equation}
	 Define
		\begin{equation*}
			T(Y)=\sum_{\substack{x,y\in Y\\x\ne y}}b_f(x,y), \qquad P(Y)=\bigl|\{e\in E(K_N):f(e)\subseteq Y\}\bigr|
		\end{equation*}
		and
		\begin{equation*}
			Z_4(Y)=\bigl|\{e\in E(K_N):e\cup f(e)\subseteq Y\}\bigr|, \quad Z_3(Y)=\bigl|\{e\in E(K_N):|(e\cup f(e))\cap Y|\ge3\}\bigr|.
		\end{equation*}
		For each edge $e=uv$, write $f(e)=xy$. Since $e\cap f(e)=\varnothing$, the edge $e$
		is counted once in each of $b_f(u,x)$, $b_f(u,y)$, $b_f(v,x)$ and $b_f(v,y)$. Thus
		\begin{equation*}
			\sum_{x\ne y}b_f(x,y)=4\binom N2<2N^2.
		\end{equation*}
		By \eqref{eq:subsetprob} and linearity of expectation, we have
		\begin{equation}\label{eq:ET}
			\mathbb ET(Y)=\sum_{x\ne y}b_f(x,y)\mathbb{P}(x,y\in Y)\le2N^2\left(\frac{2L}{N}\right)^2=8L^2.
		\end{equation}
		Similarly,
		\begin{equation}\label{eq:EP}
			\mathbb E P(Y)=\sum_{e \in E(K_N)} \mathbb{P}(f(e) \subseteq Y)\le \binom N2\left(\frac{2L}{N}\right)^2\le 2L^2.
		\end{equation}
		Since $|e\cup f(e)|=4$ for every $e\in E(K_N)$, we have
		\begin{equation}\label{eq:EZ4}
			\mathbb E Z_4(Y)=\sum_{e \in E(K_N)} \mathbb{P}(e \cup f(e) \subseteq Y)\le \binom N2\left(\frac{2L}{N}\right)^4\le \frac{8L^4}{N^2}.
		\end{equation}
		Moreover, if $|(e\cup f(e))\cap Y|\ge3$, then $Y$ contains at least one of the four three-element subsets of $e\cup f(e)$. Consequently,
		\begin{equation}\label{eq:EZ3}
			\mathbb E Z_3(Y)=\sum_{e \in E(K_N)} \mathbb{P}(|(e \cup f(e)) \cap Y| \ge 3)\le 4\binom N2\left(\frac{2L}{N}\right)^3\le \frac{16L^3}{N}.
		\end{equation}
		Applying Markov's inequality to \eqref{eq:ET}--\eqref{eq:EZ3}, we obtain
		\begin{equation*}
			\mathbb P(T(Y)>128L^2)\le\frac1{16},\qquad \mathbb P(P(Y)>32L^2)\le\frac1{16},
		\end{equation*}
		and
		\begin{equation*}
			\mathbb P\left(Z_4(Y)>\frac{128L^4}{N^2}\right)\le\frac1{16},\qquad
			\mathbb P\left(Z_3(Y)>\frac{256L^3}{N}\right)\le\frac1{16}.
		\end{equation*}
		By the union bound, there exists an $L$-element set $Y\subseteq X$ such that
		\begin{equation}\label{eq:YPZZ>}
			T(Y)\le128L^2,\quad P(Y)\le32L^2,\quad Z_4(Y)\le\dfrac{128L^4}{N^2},\quad Z_3(Y)\le\dfrac{256L^3}{N}.
		\end{equation}
	We obtain $S$ from $Y$ by the following two deletion steps (see Figure~\ref{fig:selecting-S}).
		\begin{figure}[H]
			\centering
			\begin{tikzpicture}[scale=0.95]
				
				\draw[thick] (0,0) ellipse (5.2 and 2.4);
				\node at (0,2.1) {$Y$};
				
				\draw[thick] (-3.35,0) ellipse (1.05 and 0.85);
				\node at (-3.35,1.2) {$S$};
				
				\foreach \x/\y in {
					-3.72/0.24, -3.08/0.30,
					-3.76/-0.24, -3.0/-0.30,
					-3.40/0.00
				}{
					\fill (\x,\y) circle (2pt);
				}
				
				\foreach \x/\y in {
					-1.5/1.00, -1.28/0.5, -1.62/-0.5
				}{
					\draw[thick] (\x,\y) circle (2.2pt);
				}
				
				\draw[thick] (-1.52,0.22) ellipse (0.62 and 1.10);
				\node[align=center] at (-1.52,1.68) {Step 1};
				
				\draw[dashed, thick] (1.35,0) ellipse (1.08 and 0.78);
				\node at (1.35,1.05) {$e\cup f(e)$};
				
				\fill (0.84,0.30) circle (2pt);
				\fill (1.78,0.30) circle (2pt);
				\fill (0.86,-0.30) circle (2pt);
				
				\draw[thick] (1.78,-0.30) circle (2.2pt);
				
				\draw[thick] (1.78,-0.30) ellipse (0.42 and 0.32);
				\node[align=center] at (1.35,-1.12) {Step 2};
				
				\foreach \x/\y in {
					3.10/0.92, 3.85/0.42,
					3.34/-0.45, 4.18/-0.92,
					2.82/-1.28
				}{
					\fill (\x,\y) circle (2pt);
				}
			\end{tikzpicture}
			\caption{Selecting $S$ from $Y$ after the two deletion steps.}
			\label{fig:selecting-S}
		\end{figure}
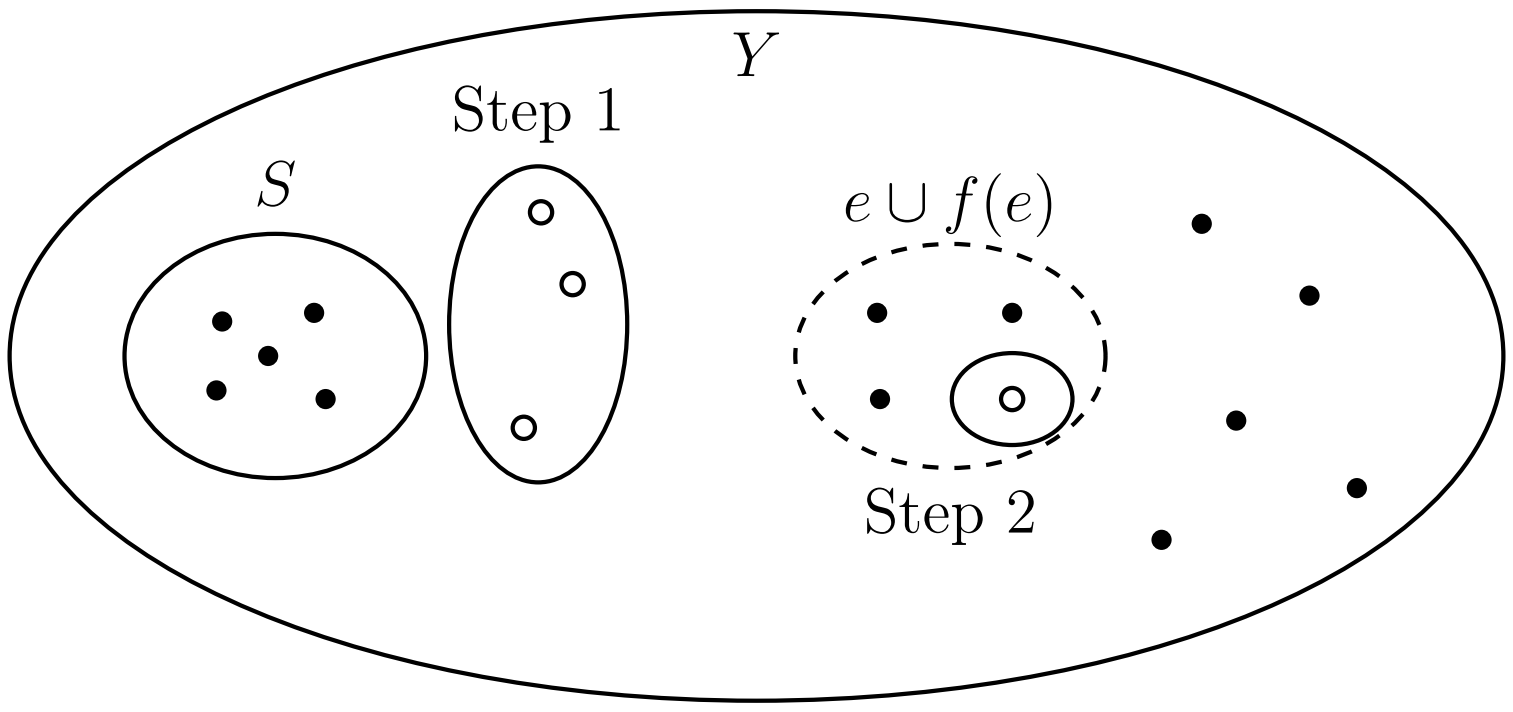
		\item  \textbf{Step 1.} For each $x\in Y$, let
		\begin{equation*}
			c_Y(x)=\sum_{y\in Y\setminus\{x\}}\bigl(b_f(x,y)+b_f(y,x)\bigr).
		\end{equation*}
		Since $\sum_{x\in Y}c_Y(x)=2T(Y)\le256L^2$, the number of vertices $x\in Y$ with
		$c_Y(x)>1024L$ is at most $\frac{2T(Y)}{1024L}\le\frac L4$. Delete all such vertices from $Y$.
		\item  \textbf{Step 2.} For each $e\in E(K_N)$ with $e\cup f(e)\subseteq Y$, delete one vertex of $e\cup f(e)$ if all four vertices remain. Since $L=4q$ and $N\ge256q^{3/2}=32L^{3/2}$, the number of additional deleted vertices is at most
		\begin{equation*}
			Z_4(Y)\le\frac{128L^4}{N^2}\le\frac{128L^4}{(32L^{3/2})^2}=\frac L8.
		\end{equation*}
		Hence, at least $L- L/4- L/8={5L}/{8}={5q}/{2}>q$ vertices remain. Choose any $q$ of these vertices to form $S$.
		
		Since $S\subseteq X$, property $\ref{c}$ holds. Since $S\subseteq Y$ and the definition of $P(Y)$ gives
		\begin{equation*}
			\bigl|\{e\in E(K_N):f(e)\subseteq S\}\bigr|\le P(Y)\le32L^2=512q^2,
		\end{equation*}
		which proves $\ref{d}$. For every $x\in S$, Step~1 gives $c_Y(x)\le1024L$. Hence
		\begin{equation*}
			\sum_{y\in S\setminus\{x\}}\bigl(b_f(x,y)+b_f(y,x)\bigr)\le c_Y(x)\le1024L=4096q,
		\end{equation*}
		which proves~$\ref{e}$. Finally, suppose that $e\subseteq S$ and $f(e)\subseteq S$ for some $e\in E(K_N)$. Then $e\cup f(e)\subseteq S\subseteq Y$ contrary to Step~2. Hence $K_N[S]$ is $f$-free, proving $\ref{a}$.
		
		To establish~\ref{b}, let
		\begin{equation*}
			W=\left\{x\in[N]\setminus S:(e\cup f(e))\setminus S=\{x\}\text{ for some }e\in E(K_N)\right\}.
		\end{equation*}
		Since $S\subseteq Y$ and $|e\cup f(e)|=4$, each vertex of $W$ is determined by an edge $e$ satisfying $|(e\cup f(e))\cap S|=3$. Hence
		\begin{equation*}
			|W|\le\bigl|\{e\in E(K_N):|(e\cup f(e))\cap S|=3\}\bigr|\le Z_3(Y)\le\frac{256L^3}{N}\le\frac N4,
		\end{equation*}
		Set $R=[N]\setminus(S\cup W)$. Since $|S|=q\le N/4$, we obtain $|R|=N-q-|W|\ge N/2$. Moreover, for every $e\in E(K_N)$, if $|(e\cup f(e))\cap S|=3$, then $(e\cup f(e))\setminus S\subseteq W$, and hence $(e\cup f(e))\cap R=\varnothing$. Thus property~$\ref{b}$ holds, completing the proof.
	\end{proof}
	We now apply Lemma~\ref{lem:host-core} to construct an $f$-free copy of $G$ and prove Theorem~\ref{thm:upper}.
	\begin{proof}[\bf Proof of Theorem~\ref{thm:upper}]
		Set $N=\left\lceil4096\left(n+m^{3/4}\right)\right\rceil$ and fix $f\in\calF_{N,0}$. Let $A=\{v\in V(G):d_G(v)>\sqrt m\}$ and $B=V(G)\setminus A$. Then $d_G(v)\le\sqrt m$ for every $v\in B$, and
		\begin{equation*}
			2m=\sum_{v\in V(G)}d_G(v)\ge\sum_{v\in A}d_G(v)>|A|\sqrt m.
		\end{equation*}
		Hence $|A|<2\sqrt m$. Let $q=\lceil2\sqrt m\rceil$. Since $m\ge1$, we have $|A|<q\le3\sqrt m$, and therefore $256q^{3/2}\le256\cdot3^{3/2}m^{3/4}<4096m^{3/4}\le N$. By Lemma~\ref{lem:host-core}, there exist disjoint sets $S,R\subseteq[N]$ satisfying $|S|=q$, $|R|\ge N/2$, and properties $\ref{a}$--$\ref{e}$.
		
		Since $|A|<|S|$, fix an injection $\phi:A\rightarrow S$. For each $v\in B$, choose $\phi(v)$ independently and uniformly from $R$. Since $\phi(A)\subseteq S$, $\phi(B)\subseteq R$ and $S\cap R=\varnothing$, and $\phi|_A$ is injective, for distinct $u,v\in V(G)$, if $\phi(u)=\phi(v)$, then $u,v\in B$.  We consider the following two types of events.
		\begin{enumerate}[label=(\roman*)]
			\item For each $\{u,v\}\in\binom{B}{2}$, define $C_{u,v}:=\{\phi(u)=\phi(v)\}$. Since $\phi(u)$ and $\phi(v)$ are chosen independently and uniformly from $R$, we have
			\begin{equation*}
				\mathbb P(C_{u,v})=\sum_{x\in R}
				\mathbb P(\phi(u)=x)\mathbb P(\phi(v)=x)=\sum_{x\in R}\frac1{|R|^2}=\frac1{|R|}.
			\end{equation*}
			\item For every ordered pair of vertex-disjoint edges $e=uv$ and $e'=xy$ of $G$, let $D_{e,e'}$ be the event that $\phi(u),\phi(v),\phi(x),\phi(y)$ are distinct and $f(\phi(u)\phi(v))=\phi(x)\phi(y)$.
		\end{enumerate}
		If none of the events $C_{u,v}$ occurs, then $\phi$ is injective. Under this assumption, if none of the events $D_{e,e'}$ occurs, then $\phi(G)$ is $f$-free, since $f(\phi(e))=\phi(e')$ implies $e\cap e'=\varnothing$.
		
		Suppose that $D_{e,e'}$ occurs. If $e\cup e'\subseteq A$, then $f(\phi(e))=\phi(e')$ contradicts property~\ref{a} of Lemma~\ref{lem:host-core}. If $|(e\cup e')\cap A|=3$, then $|(\phi(e)\cup f(\phi(e)))\cap S|=3$ and $(\phi(e)\cup f(\phi(e)))\cap R\ne\varnothing$ contrary to property~\ref{b}. Since $|e\cup e'|=4$, it follows that if $D_{e,e'}\ne\varnothing$, then $|(e\cup e')\cap B|\ge2$. Fix $v\in B$, so $d_G(v)\le\sqrt m$.	We estimate
		\begin{equation*}
		\sum_{\substack{
				(e,e'):e\cap e'=\varnothing,\,
				v\in e\cup e'}}
		\mathbb P(D_{e,e'})
		\end{equation*}
		by considering the following four cases.
		
		\item \textbf{Case 1.} Suppose that $|(e\cup e')\cap A|\le1$. If $e'\subseteq B$, then, conditional on $\phi(e)$ being an edge, the two endpoints of $e'$ are mapped independently and uniformly into $R$. Hence
		\begin{equation*}
			\mathbb P(D_{e,e'})\le\frac2{|R|^2}.
		\end{equation*}
		If $e'\cap A=\{a\}$, then $e\subseteq B$. Write $e=uv$ and $e'=aw$. Since $\phi(a)$ is fixed, the equality $f(\phi(u)\phi(v))=\phi(a)\phi(w)$ holds for at most $2a_f(\phi(a))$ ordered triples $(\phi(u),\phi(v),\phi(w))\in R^3$. Thus, by property~\ref{c} of Lemma~\ref{lem:host-core} and $|R|\ge N/2$,
		\begin{equation*}
			\mathbb P(D_{e,e'})\le\frac{2a_f(\phi(a))}{|R|^3}\le\frac{4N}{|R|^3}\le\frac8{|R|^2},
		\end{equation*}
 Moreover, for each $v\in B$, there are at most $d_G(v)m$ ordered pairs $(e,e')$ with $v\in e$,
 and at most $d_G(v)m$ with $v\in e'$. Hence $\bigl|\{(e,e'):e\cap e'=\varnothing,\ v\in e\cup e'\}\bigr|\le2d_G(v)m$. Consequently,
		\begin{equation}\label{eq:case1}
			\sum_{\substack{
					(e,e'):e\cap e'=\varnothing,\,
					v\in e\cup e'\\
					|(e\cup e')\cap A|\le1}}
			\mathbb P(D_{e,e'})
			\le \frac8{|R|^2}\cdot2d_G(v)m=
			\frac{16d_G(v)m}{|R|^2}.
		\end{equation}
		\item \textbf{Case 2.} Suppose that $|e\cap A|=2$. Since $|(e\cup e')\cap A|\le2$ whenever $D_{e,e'}\ne\varnothing$, we may assume that $e\subseteq A$ and $e'\subseteq B$. Write $e'=xy$. Since $\phi|_A$ is fixed, $f(\phi(e))$ is a fixed edge. As $\phi(x)$ and $\phi(y)$ are independent and uniformly distributed over $R$, we have $\mathbb P(D_{e,e'})\le2/{|R|^2}$. For every $v\in B$, there are at most $d_G(v)$ choices for $e'\subseteq B$ with $v\in e'$, and at most $m$ choices for $e\subseteq A$. Consequently,
		\begin{equation}\label{eq:case2}
				\sum_{\substack{
						e\subseteq A,\,e'\subseteq B,\,
						v\in e'}}
				\mathbb P(D_{e,e'})\le
				\frac{2}{|R|^2}\cdot d_G(v)m=\frac{2d_G(v)m}{|R|^2}.
		\end{equation}
		\item \textbf{Case 3.} Suppose that $|e'\cap A|=2$. Since $|(e\cup e')\cap A|\le2$ whenever $D_{e,e'}\ne\varnothing$, we may assume that $e\subseteq B$ and $e'\subseteq A$.
		Write $e=uv$. Since $\phi(e')$ is fixed, the equality $f(\phi(u)\phi(v))=\phi(e')$ holds for at most $2|\{h\in E(K_N):f(h)=\phi(e')\}|$ ordered pairs $(\phi(u),\phi(v))\in R^2$. Hence
		\begin{equation*}
			\mathbb P(D_{e,e'})\le\frac{2}{|R|^2}\bigl|\{h\in E(K_N):f(h)=\phi(e')\}\bigr|.
		\end{equation*}
		Since $\phi|_A$ is injective, the edges $\phi(e')$, for $e'\subseteq A$, are distinct.
		Thus property~\ref{d} of Lemma~\ref{lem:host-core} gives
		\begin{equation*}
			\sum_{e'\subseteq A}\bigl|\{h\in E(K_N):f(h)=\phi(e')\}\bigr|\le\bigl|\{h\in E(K_N):f(h)\subseteq S\}\bigr|\le512q^2.
		\end{equation*}
	For every $v\in B$, there are at most $d_G(v)$ choices for $e\subseteq B$ with $v\in e$. Consequently,
		\begin{align}\label{eq:case3}
			\sum_{\substack{
					(e,e'):e\cap e'=\varnothing\\
					e\subseteq B,\,e'\subseteq A,\,v\in e}}
			\mathbb P(D_{e,e'})&\le
			\frac{2d_G(v)}{|R|^2}
			\sum_{e'\in E(G[A])}
			\bigl|\{h\in E(K_N):f(h)=\phi(e')\}\bigr|\notag\\&\le
			\frac{1024d_G(v)q^2}{|R|^2}.
		\end{align}
		\item \textbf{Case 4.} Suppose that $|e\cap A|=|e'\cap A|=1$. Write $e=au$ and $e'=bw$, where $a,b\in A$ and $u,w\in B$, with $a\ne b$ and $u\ne w$. Since $\phi(a)$ and $\phi(b)$ are fixed, the equality $f(\phi(a)\phi(u))=\phi(b)\phi(w)$ holds for at most $b_f(\phi(a),\phi(b))$ ordered pairs $(\phi(u),\phi(w))\in R^2$. Hence, by \eqref{eq:bf},
		\begin{equation*}
			\mathbb P(D_{au,bw})\le\frac{b_f(\phi(a),\phi(b))}{|R|^2}.
		\end{equation*}
		Fix $v\in B$. First suppose that $v=u$. For each $b\in A$, there are at most $|N_G(b)\cap B|$ choices for $w$. By property~\ref{e} of Lemma~\ref{lem:host-core},
		\begin{align*}
			\sum_{\substack{
					a\in N_G(v)\cap A,\ 
					b\in A\setminus\{a\}\\
					w\in N_G(b)\cap B,\ w\ne v
			}}\mathbb P(D_{av,bw})&\le
			\frac1{|R|^2}
			\sum_{b\in A}|N_G(b)\cap B|
			\sum_{a\in A\setminus\{b\}}
			b_f(\phi(a),\phi(b))\\&\le
			\frac{4096q}{|R|^2}
			\sum_{b\in A}|N_G(b)\cap B|\le
			\frac{4096qm}{|R|^2},
		\end{align*}
		Similarly, when $v=w$, we have
		\begin{equation*}
	\sum_{\substack{
			b\in N_G(v)\cap A,\ a\in A\setminus\{b\}\\
			u\in N_G(a)\cap B,\ u\ne v}}
	\mathbb P(D_{au,bv})\le\frac1{|R|^2}
	\sum_{a\in A}|N_G(a)\cap B|
	\sum_{b\in A\setminus\{a\}}
	b_f(\phi(a),\phi(b))\le\frac{4096qm}{|R|^2}.
		\end{equation*}
	Consequently,
	\begin{equation}\label{eq:case4}
		\sum_{\substack{
				e,e'\in E(G),\,e\cap e'=\varnothing\\
				|e\cap A|=|e'\cap A|=1,
				v\in e\cup e'}}
		\mathbb P(D_{e,e'})
		\le\frac{8192qm}{|R|^2}.
	\end{equation}
		Combining \eqref{eq:case1}--\eqref{eq:case4}, and using $d(v)\le\sqrt m$ and $q\le3\sqrt m$, we obtain, for every $v\in B$,
		\begin{equation*}
			\sum_{\substack{(e,e'):e\cap e'=\varnothing\\v\in e\cup e'}}
			\PP(D_{e,e'})\le\frac{18d(v)m+1024d(v)q^2+8192qm}{|R|^2}
			\le\frac{60000m^{3/2}}{|R|^2}.
		\end{equation*}
		Moreover,
		\begin{equation*}
			\sum_{u\in B\setminus\{v\}}\mathbb P(C_{u,v})=\frac{|B|-1}{|R|}\le\frac n{|R|}.
		\end{equation*}
		Since $|R|\ge N/2$ and $N\ge4096(n+m^{3/4})$, it follows that
		\begin{align*}
			\sum_{u\in B\setminus\{v\}}\mathbb P(C_{u,v})+\sum_{\substack{
					(e,e'):
					e\cap e'=\varnothing\\
					v\in e\cup e'}}
			\mathbb P(D_{e,e'})&\le\frac n{|R|}+\frac{60000m^{3/2}}{|R|^2}\le\frac{2n}{N}
			+\frac{240000m^{3/2}}{N^2}\\&\le\frac2{4096}+\frac{240000}{4096^2}<\frac1{16}.
		\end{align*}
		Each event $C_{u,v}$ depends on two independent random variables. It suffices to consider the nonempty events $D_{e,e'}$. Each such event satisfies $2\le |(e\cup e')\cap B|\le4$, and hence depends on at most four and at least two random variables. By Lemma~\ref{lem:variable-lll}, we obtain
		\begin{equation*}
			\mathbb P\left(
			\bigcap_{\{u,v\}\in\binom B2}\overline{C_{u,v}}\cap
			\bigcap_{\substack{
					(e,e')\in E(G)^2\\
					e\cap e'=\varnothing}}
			\overline{D_{e,e'}}
			\right)>0.
		\end{equation*}
		Consequently, there exists an injective mapping $\phi:V(G)\to[N]$ such that $f(\phi(e))\notin E(\phi(G))$ for every $e\in E(G)$. Hence $\phi(G)$ is $f$-free, completing the proof.
	\end{proof}
	\section{The Lower Bound}\label{sec:The Lower Bound}
	In this section, we prove Theorem~\ref{thm:lower-all} by constructing graphs with $n$ vertices and $m$ edges. We distinguish two cases according to whether $m<2n$ or $m\ge2n$.
	\begin{proof}[\bf Proof of Theorem~\ref{thm:lower-all}]
		Let $n\ge2$ and $\lceil n/2\rceil\le m\le\binom n2$. For $2\le n\le7$, choose any graph $G_*$ with $n$ vertices, $m$ edges and no isolated vertices. Since
		\begin{equation*}
	m^{3/4}\le\left(\frac{n^2}{2}\right)^{3/4}\le\frac{\sqrt7}{2^{3/4}}n<2n,
		\end{equation*}
		we have 
		\begin{equation*}
			g_0(G_*)\ge n\ge\frac13(n+m^{3/4}).
		\end{equation*}
		Hence we may assume that $n\ge8$.
		\item \textbf{Case 1.} $m<2n$. If $n$ is even, take a perfect matching on $n$ vertices. If $n$ is odd, take $P_3\,\dot\cup\,\frac{n-3}{2}K_2$. In either case, the graph has no isolated vertices and exactly $\lceil n/2\rceil$ edges. Add edges to obtain a graph $G_*$ with exactly $m$ edges. Since $n\ge8$ and $m<2n$, we have $m^{3/4}<(2n)^{3/4}\le n$. Consequently, 
		\begin{equation*}
			g_0(G_*)\ge n\ge\frac12(n+m^{3/4}).
		\end{equation*}
		\item \textbf{Case 2.} $m\ge2n$. Let $s=\lfloor\sqrt m\rfloor$. Since $16\le m\le\binom n2$, we have $4\le s<n$ and $s\ge{\sqrt m}/{2}$. Start with $K_s$ and join each of the remaining $n-s$ vertices to a fixed vertex of $K_s$. The resulting graph has no isolated vertices, and at most
		\begin{equation*}
			\binom s2+n-s\le\frac m2+n\le m.
		\end{equation*}
		edges. Add edges to obtain a graph $G_*$ with exactly $m$ edges and $K_s\subseteq G_*$. By Proposition~\ref{prop:mono} and Corollary~\ref{cor:clique}, we obtain
		\begin{equation*}
			g_0(G_*)\ge g_0(K_s)\ge\frac{s^{3/2}}{200}\ge\frac{m^{3/4}}{200\cdot2^{3/2}}>\frac{m^{3/4}}{600}.
		\end{equation*}
		Together with $g_0(G_*)\ge n$, this gives
		\begin{equation*}
			g_0(G_*)\ge\max\left\{n,\frac{m^{3/4}}{600}\right\}\ge\frac1{1200}(n+m^{3/4}).
		\end{equation*}
		The proof is complete.
	\end{proof}
	Combining Theorems~\ref{thm:upper} and~\ref{thm:lower-all}, we obtain the following result for graphs with given numbers of vertices and edges.
	\begin{corollary}\label{cor:final}
		For every $n\ge2$ and every integer $m$ satisfying $\lceil n/2\rceil\le m\le\binom n2$, there exists a graph $G_*$ with $n$ vertices, $m$ edges and no isolated vertices such that
		\begin{equation*}
			g_0(G_*)=\Theta\bigl(n+m^{3/4}\bigr).
		\end{equation*}
	\end{corollary}
	\begin{proof}[\bf Proof]
		The upper bound follows from Theorem~\ref{thm:upper}, and the lower bound follows from Theorem~\ref{thm:lower-all}.
	\end{proof}
	
\section*{Declaration on the use of AI}
	The author used ChatGPT to assist in discussing proof strategies, checking proofs, and improving exposition.


\begin{thebibliography}{99}
	\bibitem{AlonCaro1986}
	N. Alon and Y. Caro,
	\newblock Extremal problems concerning transformations of the
	set of edges of the complete graph,
	\newblock \emph{European J. Combin.} \textbf{7} (1986) 93--104.
	
	\bibitem{AlonCaroTuza1989}
	N. Alon, Y. Caro, and Z. Tuza,
	\newblock Sub-Ramsey numbers for arithmetic progressions,
	\newblock \emph{Graphs Comb.} \textbf{5} (1989) 307--314.
	
	\bibitem{AlonSpencer2016}
	N. Alon and J. H. Spencer,
	\newblock \emph{The Probabilistic Method},
	\newblock 4th ed., Wiley, 2016.
	
	
	\bibitem{Caro1987}
	Y. Caro,
	\newblock Extremal problems concerning transformations of
	the edges of the complete hypergraphs,
	\newblock \emph{J. Graph Theory} \textbf{11} (1987) 25--37.
	
\bibitem{CPTV2026}
Y. Caro, B. Patk\'os, Z. Tuza, and M. Vizer,
\newblock Edge mappings of graphs: Ramsey-type parameters,
\newblock arXiv preprint arXiv:2402.01004v2, 2026.
	
	\bibitem{CPTV2025}
	Y. Caro, B. Patk\'os, Z. Tuza, and M. Vizer,
	\newblock Edge mappings of graphs: Tur\'an type parameters,
	\newblock \emph{European J. Combin.} \textbf{127} (2025) 104140.
	
	\bibitem{CaroSchonheim1988}
	Y. Caro and J. Sch\"onheim,
	\newblock On edge-mappings with fixed edges avoiding free triangles,
	\newblock \emph{Ars Combin.} \textbf{25-A} (1988) 159--164.
	
	
	\bibitem{CFS2016}
	D. Conlon, J. Fox, and B. Sudakov,
	\newblock Short proofs of some extremal results II,
	\newblock \emph{J. Combin. Theory Ser. B}
	\textbf{121} (2016) 173--196.
	

	
	\bibitem{ErdosHajnal1958}
	P. Erd\H{o}s and A. Hajnal,
	\newblock On the structure of set-mappings,
	\newblock \emph{Acta Math. Acad. Sci. Hungar.}
	\textbf{9} (1958) 111--131.
	
	
	\bibitem{Furedi1991}
	Z. F\"uredi,
	\newblock Maximal independent subsets in Steiner systems
	and in planar sets,
	\newblock \emph{SIAM J. Discrete Math.}
	\textbf{4} (1991) 196--199.
	
	\bibitem{GJS2026}
	L. Gishboliner, Z. Jin, and B. Sudakov,
	\newblock Set mappings for general graphs,
	\newblock arXiv preprint arXiv:2601.00766, 2026.
	
\bibitem{Spencer1972}
J. Spencer,
\newblock Tur\'an's theorem for $k$-graphs,
\newblock \emph{Discrete Math.} \textbf{2} (1972), 183--186.
	\end{thebibliography}
\end{document}